\documentclass{ifacconf}

\usepackage{graphicx}      
\usepackage{natbib}        
\usepackage{latexsym,amsmath,amssymb,amsopn, graphicx}
\usepackage{enumerate}
\usepackage{array} 
\usepackage{float} 
\usepackage{multirow} 
\usepackage{multicol} 
\usepackage{array}
\usepackage[normalsize,it,hang]{subfigure}
\usepackage{tabularx}
\begin{document}
\begin{frontmatter}

\title{Reweighted nuclear norm regularization: A SPARSEVA approach\thanksref{footnoteinfo}} 

\thanks[footnoteinfo]{This work was partially supported by the Swedish Research Council and the Linnaeus Center ACCESS at KTH. }

\author[First]{Huong Ha} 
\author[First]{James S. Welsh} 
\author[Second]{Niclas Blomberg}
\author[Second]{Cristian R. Rojas}
\author[Second]{Bo Wahlberg}

\address[First]{School of Electrical Engineering and Computer Science, The University of Newcastle, Australia (e-mail: huong.ha@uon.edu.au, james.welsh@newcastle.edu.au).}
\address[Second]{Department of Automatic Control and ACCESS, School of Electrical Engineering, KTH Royal Institute of Technology, SE-100 44 Stockholm, Sweden.(e-mail: nibl@kth.se, cristian.rojas@ee.kth.se, bo.wahlberg@ee.kth.se).}

\begin{abstract}                
The aim of this paper is to develop a method to estimate high order FIR and ARX models using least squares with re-weighted nuclear norm regularization. Typically, the choice of the tuning parameter in the reweighting scheme is computationally expensive, hence we propose the use of the SPARSEVA (SPARSe Estimation based on a VAlidation criterion) framework to overcome this problem. Furthermore, we suggest the use of the prediction error criterion (PEC) to select the tuning parameter in the SPARSEVA algorithm. Numerical examples demonstrate the veracity of this method which has close ties with the traditional technique of cross validation, but using much less computations.
\end{abstract}


\end{frontmatter}

\section{Introduction}
\vspace{-0.1cm}
In the last few decades, utilizing regularization in estimating transfer functions has received much attention (\cite{Pillonetto:14}). The basic idea is to add a penalty on the weighted size of the parameters to the least square cost function. There are several different regularization schemes that have been examined. For example, Tikhonov regularization (or ridge regression) uses a weighted $l_2$ norm as the penalty function to solve ill-posed problems in \cite{Tikhonov:77}. More recently, an important contribution to this field is from \cite{Pillonetto:10}, where estimation methods based on Gaussian processes and spline kernels were examined. This approach has shown to have close ties with $l_2$ regularization in \cite{Chen:11}. Another purpose of regularization is to obtain a sparse estimate. One of the early contributions to this approach is the LASSO (Least Absolute Shrinkage and Selection Operation) that minimizes the least square cost function under a constraint on the $l_1$ norm of the parameter vector, which could also be interpreted as $l_1$ regularization in \cite{Tibshirani:96}. The idea behind of the LASSO is to use the $l_1$ norm as a surrogate for the $l_0$ norm.

Recently, the nuclear norm has been used in this approach and shown to be a useful method for identifying low order linear models. It is first suggested to be used in sparse estimation as a rank minimization heuristic in \cite{Fazel:01}. After that, \cite{Grosmann:09a, Grosmann:09b} have proposed a system identification method to fit a flexible FIR model to measured data and at the same time penalize model order complexity by minimizing the nuclear norm of the corresponding Hankel matrix. An advantage of this approach compared to classical system identification methods is that the corresponding optimization problem is convex. The trade-off between fit and complexity is determined by a design parameter, which is typically found by cross-validation. This has been extended to ARX models by \cite{Hakan:12}.


As the nuclear norm is an approximation of the rank function, by using a singular value weighting, this approximation can be further improved (\cite{Fazel:03},\cite{Mohan:10}). The objective of this paper is to study the use of re-weighting for the nuclear norm system identification problem by developing a method to estimate high order FIR and ARX models using least squares with re-weighted nuclear norm regularization. 

The choice of the tuning parameter will then, however, be more involved as the reweighted nuclear norm regularization is an iterative process. To overcome this problem we propose the use of the SPARSEVA framework, which eliminates  the  need  for cross-validation  to  tune  the  regularization  parameter (\cite{Rojas:11, Rojas:14}).
The SPARSEVA  is originally motivated by classical structure selection criteria, such as Akaike Information Criteria (AIC) and Bayesian Information Criteria (BIC). Here, we propose the use of the prediction error criterion, which leads to the SPARSEVA - PEC method.

The paper is organized as follows. Section 2 describes the system identification problem. Section 3 establishes the proposed reweighted nuclear norm regularization method and discusses the use of SPARSEVA. Section 4 suggests an alternative way to choose the regularization parameter in the SPARSEVA method. In section 5, some numerical experiments and results are presented. Finally, section 6 provides the conclusions.

\vspace{0.01cm}
\section{System Identification Problem}
\vspace{-0.1cm}

Consider a discrete-time linear, time invariant, single input single output stable system,
\begin{equation}\label{eq-System}
	y(t) = G_0(q)u(t) + H_0(q)e(t), 
\end{equation}
where $u(t)$ is the input signal and $e(t)$ is zero mean white noise, independent of the input.

Given past data $\{y(t),u(t)\}$, $ t=1,...,N$ the task is to estimate a model of the transfer function $G_o(q)$. A typical approach is to apply the Prediction Error Method (PEM), however, PEM can have problems with local minima, e.g. when a rational transfer function model is estimated using the Output Error (OE) model. The advantage of the PEM is that only a few parameters need to be estimated.

Another common and simple approach to estimate $G_0(q)$ in (\ref{eq-System}) is to use the Finite Impulse Response (FIR) model,
\begin{equation}
G(q) = \sum\limits_{k=1}^{n} g_kq^{-k}, 
\end{equation}
where $n$ is the order of the FIR model. The vector $\theta = [g_1, g_2, ..., g_{n}]^T$ can then be estimated using least squares by solving the linear regression problem corresponding to the model,
\begin{equation} \label{eq_regression}
y(t) = \varphi(t)^T\theta,
\end{equation}
with
\begin{equation}
\varphi(t) = [u(t-1),\ ..., u(t-n)]^T.
\end{equation}

However, as $n$ increases, the variance of the estimate will also increase. Regularization can then be used to overcome this issue (\cite{Ljung:91})(\cite{Chen:11}). As mentioned in the introduction, there are many regularization schemes where in this paper, we will only focus on the nuclear norm regularization.

Consider the Hankel matrix $\mathcal{H}(\theta)$ formulated from $\theta$,
\begin{equation} \label{Hankel_mt}
\mathcal{H}(\theta) =  \begin{bmatrix}
    g_1 & g_2 & \dots  & g_m \\
    g_2 & g_3 & \dots  & g_{m+1} \\
    \vdots & \vdots & \ddots & \vdots \\
    g_m & g_{m+1} & \dots  & g_n
\end{bmatrix}
\end{equation}
where $m=(n+1)/2$ (n should be chosen as an odd number).
For a linear system with impulse responses $\theta = [g_1, g_2, ..., g_{n}]^T$, the rank of the Hankel matrix in (\ref{Hankel_mt}) equals to the order of the linear system (\cite{Kailath:80}). Therefore, the Hankel matrix rank minimization is often used as a important approach in model order reduction problems, of which a common surrogate heuristic for the rank function is the nuclear norm.

The nuclear norm for a matrix $X \in \mathbb{R}^{n\times m}$ is defined by
\begin{equation}
\parallel X \parallel_* = \sum\limits_{i=1}^{min(n,m)} \sigma_i(X),
\end{equation}
where ${\sigma_i(X)}$ are the singular values of $X$. Note that the nuclear norm gives a measure of the matrix rank and is equal to this rank if the singular values equal to 1 and 0.

The nuclear norm regularization algorithm can then be written in the following form by semi-definite program (SDP) (\cite{Grosmann:09b}):
\begin{equation} \label{eq-FIR_nuc}
\begin{aligned}
& \underset{\theta, W_1, W_2}{\text{min}} 
& & \sum\limits_{t=1}^{N} \Big({y(t) - \varphi(t)^T\theta}\Big)^2 +\ \dfrac{\lambda}{2} \text{Tr}(W_1+W_2) \\
& \ \ \ \text{s.t.}
& & \ \begin{bmatrix}
  W_1 \ \ \ \ \ \ \ \ \ \mathcal{H}(\theta) \\
  \mathcal{H}^T(\theta) \ \ \ \ \ \ \ W_2
  \end{bmatrix} \geq 0. \\
\end{aligned}
\end{equation}

The regularization parameter $\lambda$ can be chosen using cross validation, where the data is split into two parts, namely the estimation part and the validation part. The model is then estimated from the estimation data using different values of $\lambda$. The estimates are then validated using the validation data; the value of $\lambda$ that gives the smallest sum square error will be chosen. The model can then be re-estimated using the whole dataset with this chosen value of $\lambda$. 

In addition, another extension of the nuclear norm regularization (\cite{Hakan:12}) is to use the ARX model,
\begin{equation} \label{ARX}
	y(t) = \dfrac{B(q)}{A(q)}u(t) + \dfrac{1}{A(q)}e(t), 
\end{equation}
with 
\begin{equation} \nonumber
\begin{aligned}
& B(q) = b_1q^{-1}+...+b_{n_B}q^{-n_B}, \\
& A(q) = 1+a_1q^{-1}+...+a_{n_A}q^{-n_A}.
\end{aligned}
\end{equation}
where $n_A$ is the order of the ARX model.

As shown by \cite{Ljung:92}, when $n_A$ increases, $G_0(q)$ can be well approximated by $B(q)/A(q)$ and $H_0(q)$ can be approximated by $1/A(q)$. As demonstrated in \cite{Hakan:12}, using a nuclear norm relaxation can allow the high order estimated model to have a similar accuracy to a low order model; and the performance of this method is competitive with respect to the prediction error method.

For an ARX model, the linear regression problem to estimate the system parameters can be formulated as
\begin{equation} \label{eq_regressionARX}
A(q)y(t) = B(q)u(t)+e(t),
\end{equation}
which is the same as
\begin{equation}
y(t) = -\sum\limits_{k=1}^{n_{A}} a_ky(t-k) + \sum\limits_{k=1}^{n_{B}} b_ku(t-k)+e(t).
\end{equation}

By using nuclear norm regularization, the algorithm to estimate the parameter vector $\theta$ is as follows:
\begin{equation} \label{eq-ARX_nuc}
\begin{array}{cl}
\begin{aligned}
& \underset{\theta, W_{1a}, W_{2a},W_{1b}, W_{2b}}{\text{min}} \ \ \ \ \sum\limits_{k=1}^{N}\Big(y(t)-\varphi^T\theta\Big)^2 \\
& \ \ \ \ \ \ \ \ \ \ \ \ \ \ \ \ +\ \dfrac{\lambda_a}{2} \text{Tr}(W_{1a}+W_{2a})+\ \dfrac{\lambda_b}{2} \text{Tr}(W_{1b}+W_{2b})\\
& \ \ \ \ \ \ \text{s.t.} \ \ \ \ \ \ \ \ \ \ \begin{bmatrix}
  W_{1a} \ \ \ \ \ \ \ \ \ \mathcal{H}(a) \\
  \mathcal{H}^T(a) \ \ \ \ \ \ \ W_{2a}
  \end{bmatrix} \geq 0, \\
& \ \ \ \ \ \ \ \ \ \ \ \ \ \ \ \ \ \ \ \ \begin{bmatrix}
  W_{1b} \ \ \ \ \ \ \ \ \ \mathcal{H}(b) \\
  \mathcal{H}^T(b) \ \ \ \ \ \ \ W_{2b}
  \end{bmatrix} \geq 0,
\end{aligned}
\end{array}
\end{equation}
where \begin{equation} \nonumber
\begin{aligned}
& \varphi(t) = [-y(t-1)\ ...\ -y(t-{n_A})\ u(t-1)\ ...\ u(t-{n_B})]^T, \\
& \theta \ \ \ = [a_1\ \ ...\ a_{n_A}\ b_1\ ...\ b_{n_B}]^T, \\
& a = [a_1\ \ ...\ a_{n_A}]^T,\ \ b = [b_1\ ...\ b_{n_B}]^T.
\end{aligned}
\end{equation}

Another contribution to the nuclear norm regularization was recently suggested by \cite{Blomberg:14, Blomberg:15}. Here, a new methodology that utilizes the regularization path to approximately calculate all solutions of the nuclear norm regularization optimization problem to a certain tolerance of the model reduction problem as a function of the design parameter. 

\vspace{0.2cm}
\section{Reweighted nuclear norm regularization with SPARSEVA}
\vspace{-0.1cm}

The reweighted nuclear norm was developed by \cite{Mohan:10} and \cite{Fazel:03} for solving the rank minimization problem. The idea is to use the log-det heuristic to approximate the rank function. Based on this, we develop a regularization method that utilizes the reweighted nuclear norm to make a further improvement on the nuclear norm regularization.

Consider the linear regression problem,
\begin{equation} \label{lin-Reg}
Y_N = \Phi_N^T\theta + e,
\end{equation}
where $\theta^{n\times1}$ is the parameter vector to be estimated.

The regularization method using the reweighted nuclear norm for the FIR model can be written as (\cite{Mohan:10}),
\begin{equation} \label{eq-FIR_logdet}
\begin{aligned}
& \underset{\theta, W_1, W_2}{\text{min}} \ \ \ \ \parallel Y_N - \Phi_N^T\theta \parallel_2^2 \ +\ \dfrac{\lambda}{2}\ \text{logdet}(W_1+\delta I) \\
& \ \ \ \ \ \ \ \ \ \ \ \ \ \ \ +\ \dfrac{\lambda}{2}\ \text{logdet}(W_2+\delta I) \\
& \ \ \ \text{s.t.} \ \ \ \ \ \ \begin{bmatrix}
  W_1 \ \ \ \ \ \ \ \ \ \mathcal{H}(\theta) \\
  \mathcal{H}^T(\theta) \ \ \ \ \ \ \ W_2
  \end{bmatrix} \geq 0, \\
\end{aligned}
\end{equation}
where $\delta$ is a small additional regularization constant. And $\lambda$ is the regularized parameter. 

The previous convex optimization can be solved locally using an iterative process, where the $k^{th}$ step is to solve the following problem (\cite{Mohan:10}),
\begin{equation} \label{eq-logdetIter1}
\begin{aligned}
& \underset{\theta, W_1, W_2}{\text{min}} \ \parallel Y_N - \Phi_N^T\theta \parallel_2^2 \ +\ \dfrac{\lambda^k}{2} \text{Tr}(W_1^k+\delta I)^{-1}W_1\\
& \ \ \ \ \ \ \ \ \ \ \ \ \ \ \  +\ \dfrac{\lambda^k}{2} \text{Tr}(W_2^k+\delta I)^{-1}W_2 \\
& \ \ \ \text{s.t.} \ \ \ \ \begin{bmatrix}
  W_1 \ \ \ \ \ \ \ \ \ \mathcal{H}(\theta) \\
  \mathcal{H}^T(\theta) \ \ \ \ \ \ \ W_2
  \end{bmatrix} \geq 0 \\
\end{aligned}
\end{equation}

However, it is non-trivial to update $\lambda$ using cross validation technique in every iteration. One way to avoid this problem is to use the SPARSEVA method as developed by \cite{Rojas:11, Rojas:14}. This method essentially removes the use of the regularization parameter $\lambda$ by introducing a constant $\epsilon_N$ which is based on some validation criterion. Specifically, SPARSEVA replaces the traditional $l_p$ norm regularization criterion
\begin{equation}
\begin{aligned}
& \underset{\theta}{\text{min}}
& & \parallel Y_N - \Phi_N^T\theta \parallel _2^2 \\
& \ \text{s.t.}
& & \parallel\theta\parallel_p\ \leq \eta
\end{aligned}
\end{equation}
by the criterion 
\begin{equation} \label{eq-VN}
\begin{aligned}
& \underset{\theta}{\text{min}}
& & \parallel\theta\parallel_p \\
& \ \text{s.t.}
& & \ V_N(\theta) \leq V_N(\hat{\theta}_{LS})(1+\epsilon_N)
\end{aligned}
\end{equation}
where $\parallel . \parallel_p$ is the $l_p$ norm, $V_N(\theta) = \parallel Y_N - \Phi_N^T\theta \parallel _2^2$, $\hat{\theta}_{LS} = (\Phi_N\Phi_N^T)^{-1}\Phi_NY_N$ and $\epsilon_N > 0$.  

As suggested in the existing SPARSEVA method, $\epsilon_N$ can be chosen as $2n/N$ or $\text{log}(N)n/N$ based on AIC and BIC respectively (\cite{Rojas:14}).

Finally, by using the SPARSEVA method, the implementation of the reweighted nuclear norm can be written as,
\begin{equation} \label{eq-logdetIter2}
\begin{aligned}
& \underset{\theta, W_1, W_2}{\text{min}} \ \ \ \ \text{Tr}(W_1^k+\delta I)^{-1}W_1 +\ \text{Tr}(W_2^k+\delta I)^{-1}W_2 \\
& \ \ \ \text{s.t.} \ \ \ \ \ \ \begin{bmatrix}
  W_1 \ \ \ \ \ \ \ \ \ \mathcal{H}(\theta) \\
  \mathcal{H}^T(\theta) \ \ \ \ \ \ \ W_2
  \end{bmatrix} \geq 0 \\
& \ \ \ \ \ \ \ \ \ \ \ \ \ \ V_N(\theta) \leq V_N(1+\epsilon_N).
\end{aligned}
\end{equation}

\vspace{0.1cm}
\section{The SPARSEVA-PEC}
\vspace{-0.1cm}

For the SPARSEVA method, the key question is how to choose the regularization parameter $\epsilon_N$. In this section, we will describe an alternative to select the value of $\epsilon_N$ as to that proposed in \cite{Rojas:14}.

In general, we want to choose $\epsilon_N$ so that the true parameter vector of the linear regression (\ref{lin-Reg}), $\theta_0$, belongs to the set $ V_N(\theta_o)\leq V_N(\hat{\theta}_{LS})(1+\epsilon_N)$ with a high probability. We also need this set to be as small as possible in order to guarantee a good fit to the data, e.g. $\epsilon_N$ needs to be small.

We use the results from Sec. 11.5 in \cite{Soderstrom:89} to choose $\epsilon_N$. From Eq. (11.53) in \cite{Soderstrom:89}, we have:
\begin{equation} \label{V0_neweN}
V_N(\theta_o)\approx \frac{V_N(\hat{\theta}_{LS})}{(1-n/N)}.
\end{equation}
The right hand side of Eq. (\ref{V0_neweN}) is equal to an unbiased estimate of the true prediction error variance $\sigma_e^2$. Eq. (\ref{V0_neweN}) suggests a way to choose $\epsilon_N$, that is to set
$$
V_N(\hat{\theta}_{LS})(1+\epsilon_N) = \frac{V_N(\hat{\theta}_{LS})}{(1-n/N)},
$$
which means,
\begin{equation} \label{eN_new}
\epsilon_N=\frac{1}{(1-n/N)}\frac{n}{N}.
\end{equation}
Note that for small $n/N$, $ \epsilon_N \approx \frac{n}{N}$. We will denote this choice of $\epsilon_N$ by the PEC. In the numerical demonstration of this paper, we use the exact value $1/(1-n/N)\times n/N$ for $\epsilon_N$.

Moreover, if we instead evaluate $\hat{\theta}_{LS}$ on future fresh data, the parsimony principle leads to the final prediction error (FPE) criterion,
$$
V_N(\theta_o)\approx V_N(\hat{\theta}_{LS})\frac{(1+n/N)}{(1-n/N)},
$$
as mentioned in Eq. (11.54) in \cite{Soderstrom:89}. Hence, in this case,
\begin{equation}
\epsilon_N=\frac{1}{(1-n/N)}\frac{2n}{N}\approx \frac{2n}{N},
\end{equation}
which corresponds to the AIC. We can see that the SPARSEVA - PEC gives the smallest set while the SPARSEVA - AIC set is slightly larger and the SPARSEVA - BIC ($\epsilon_N=\text{log}(N)n/N$) set is the largest one. 

A further understanding of these design choices can be interpreted as following. From \cite{Soderstrom:89}, it is shown that
 $$
\frac{N[V_N(\theta_o)-V_N(\hat{\theta}_{LS})]}{[V_N(\hat{\theta}_{LS})/(1-n/N)]}
$$ is $\mathcal{X}_2 $ distributed with $n$ degrees of freedom. Therefore, in general, given a probability $\alpha$, there exists a value of $\epsilon_N$ so that the probability that $\theta_0$ belongs to the set $V_N(\theta_o)\leq V_N(\hat{\theta}_{LS})(1+\epsilon_N)$ is $\alpha$.

Note that this set is indeed a confidence set for $\theta_0$, and could be obtained from the result that $\sqrt{N}[\hat{\theta}_{LS}-\theta_o]$ is Gaussian distributed with zero mean and co-variance matrix $\sigma_e^2[\frac{1}{N}\Phi_N \Phi^T_N]^{-1}$. For the FIR case with deterministic regression vector, this is exact, but for ARX it is asymptotic. Note again that $V_N(\hat{\theta}_{LS})/(1-n/N)$, which corresponds to the choice of PEC, is an unbiased estimate of the variance $\sigma_e^2$ when $e$ is white noise.

In summary, we can see that the BIC is asymptotically consistent while the AIC typically gives a too small set. The PEC gives an even smaller set that gives more trust to the least squares estimate. However, particularly when we are not in the asymptotic (large N) regime, the PEC is a good choice as confirmed by the numerical example in Sec. 5.

\vspace{0.2cm}
\section{Numerical results}
\vspace{-0.1cm}

In this section, numerical examples are presented to demonstrate the effectiveness of the proposed method. In Sec. 5.1, we will consider the case when the output disturbance is white noise whilst in Sec. 5.2, a coloured noise output disturbance will be considered.

\vspace{-0.1cm}
\subsection{Output disturbance is white noise}
\vspace{-0.1cm}

We conduct an experiment similar to the experiment in \cite{Hakan:12}, i.e., 150 systems with order varying from 1 to 10 are generated using the command $\tt{drss}$ from Matlab. All the systems will have poles with magnitude less than 0.9. White noise is added to the system output at 4 different noise levels corresponding to the Cram\'{e}r-Rao bound of the model fit of 90\%, 77\%, 68\% and 55\%. For each noise level, the systems are estimated with three different input excitation signal and output noise realizations. The input signals have a low-pass characteristic, and is generated by filtering zero mean Gaussian white noise with unit variance through the filter
\begin{equation} \nonumber
F_u(q) = \dfrac{0.436}{1-0.9q^{-1}}.
\end{equation}
The sample size in each run is $N=450$ and the number of parameters $n$ of the FIR model is 35. 

The model fit is recorded for each system and then compared between several methods. The model fit is defined as
\begin{equation} \label{fit}
W = 100 \left(1 - \left[\dfrac{\sum_{k=1}^{n}(g_{k}-\widehat{g}_{k})^{2}}{\sum_{k=1}^{n}(g_{k}-\text{mean}(g_{k}))^{2}}\right]^{1/2}\right)
\end{equation}
where $ g_{k} $ is the impulse response of the real system and $\widehat{g}_{k}$ is the impulse response of the estimated system and $\text{mean}(g_k)=1/n\sum_{k=1}^{n}g_{k}$. 

The four methods that we compare are: 
\begin{enumerate}
\item
The nuclear norm regularization for an FIR model with cross validation (CV-FIR-N).
\item
The nuclear norm regularization for an FIR model with SPARSEVA-PEC (SPe-FIR-N).
\item
The reweighted nuclear norm regularization for an FIR model with SPARSEVA-PEC (SPe-FIR-RN).
\item
The prediction error method (PEM) that uses the output-error model where it is assumed that the order of the system is already known. The PEM estimates are generated by the command \texttt{oe} from the $\text{MATLAB}^{\text{TM}}$ System Identification Toolbox (R2014b).
\end{enumerate}

\begin{figure}[H]
\begin{center}
\includegraphics[width=57mm,height=36mm]{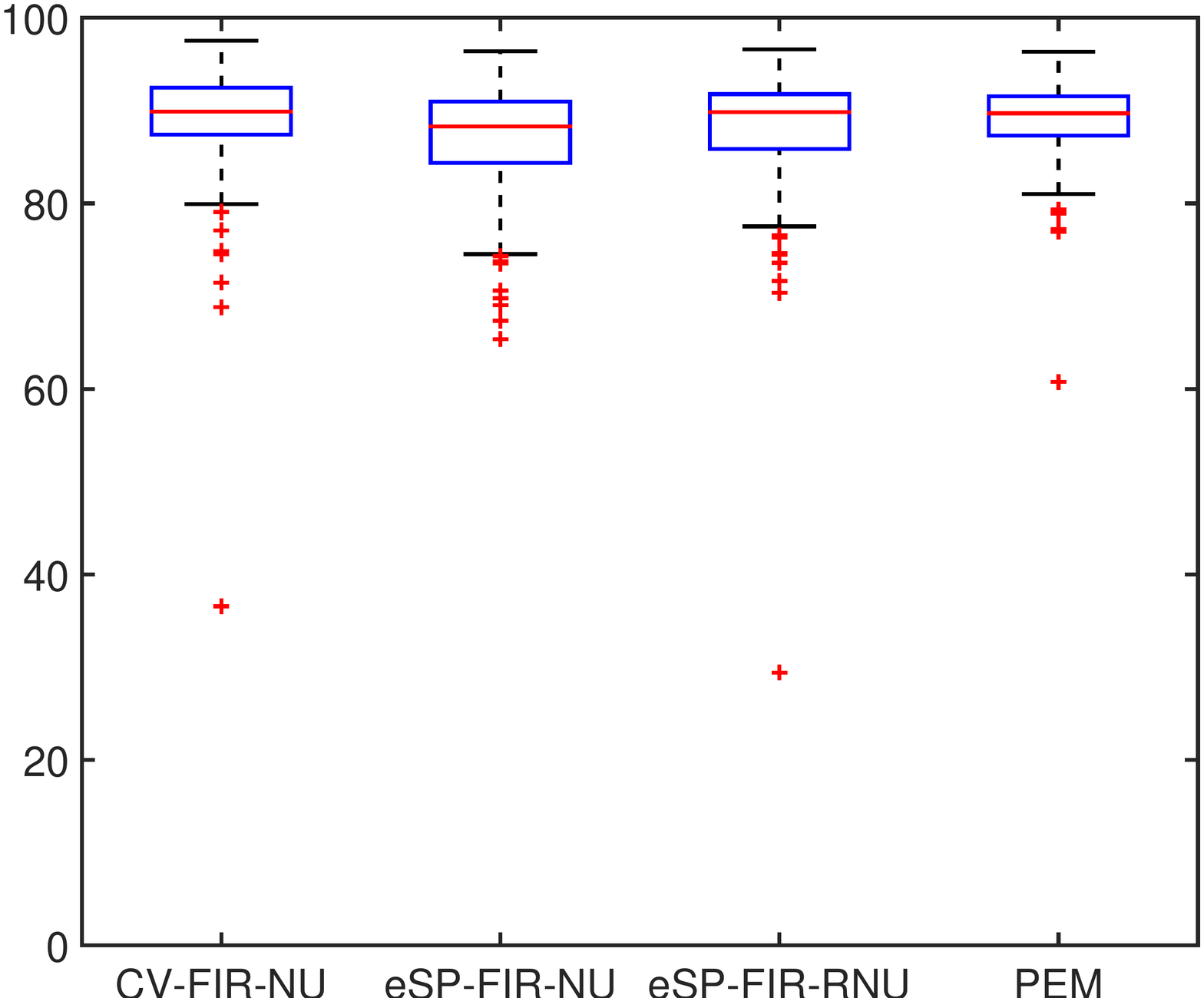}
\caption{Boxplot for Cram\'{e}r-Rao bound of 90\%.} 
\label{Boxplotl1} 
\end{center}
\end{figure}

\begin{figure}[h]
\begin{center}
\includegraphics[width=57mm,height=36mm]{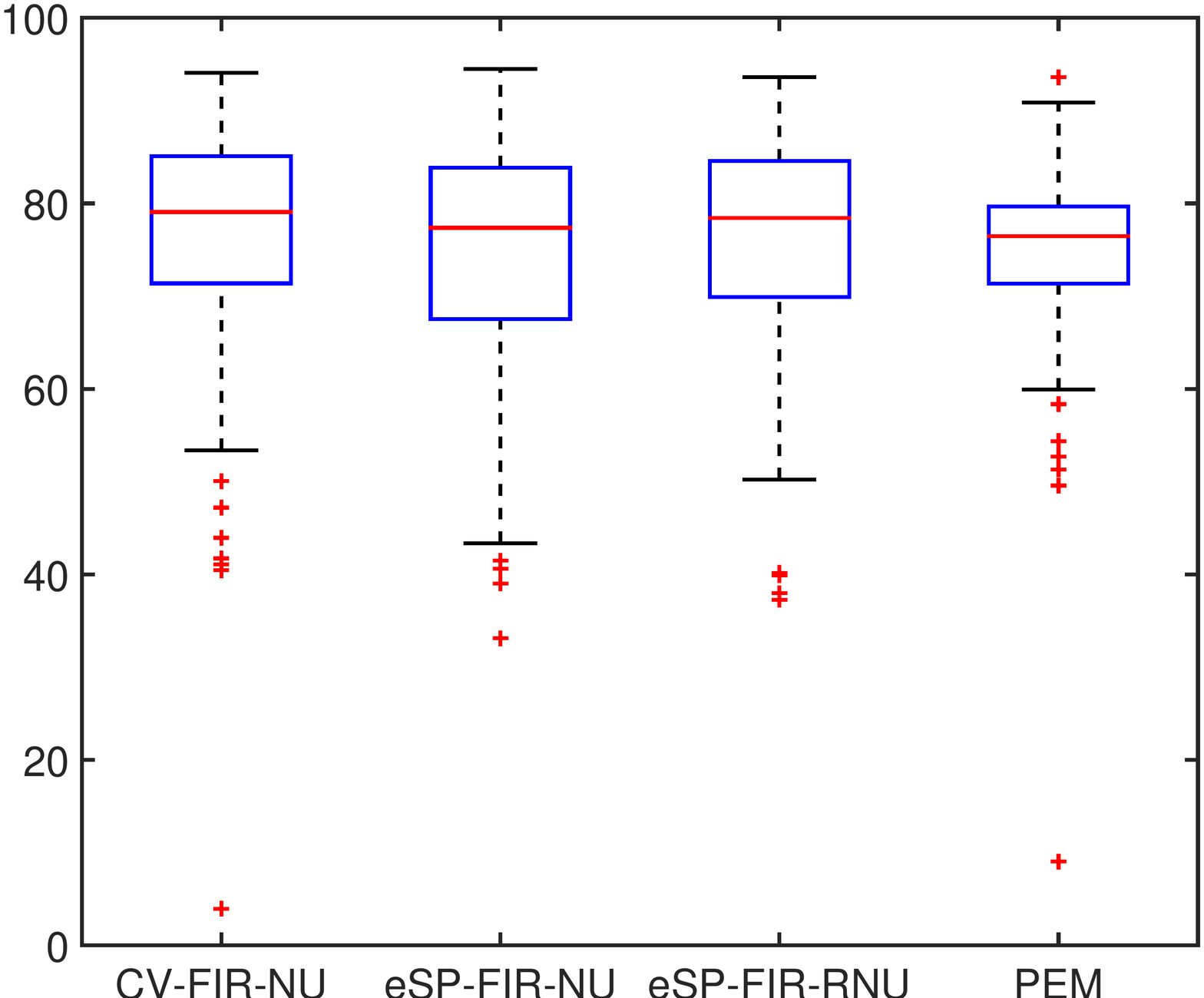}
\caption{Boxplot for Cram\'{e}r-Rao bound of 77\%.} 
\label{Boxplotl2} 
\end{center}
\end{figure}

\begin{figure}[h]
\begin{center}
\includegraphics[width=57mm,height=36mm]{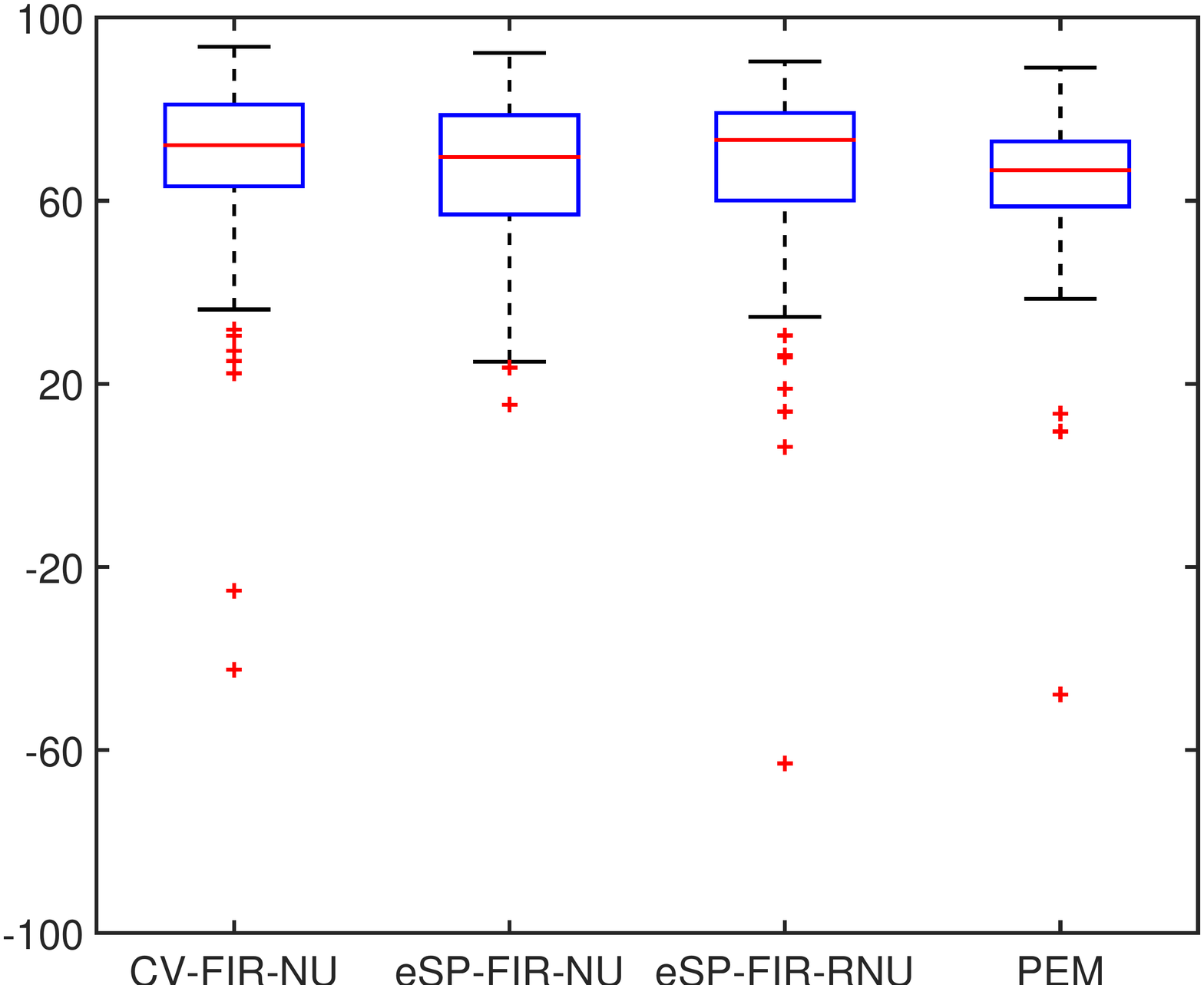}
\caption{Boxplot for Cram\'{e}r-Rao bound of 68\%.} 
\label{Boxplotl3} 
\end{center}
\end{figure}

\begin{figure}[h]
\begin{center}
\includegraphics[width=57mm,height=36mm]{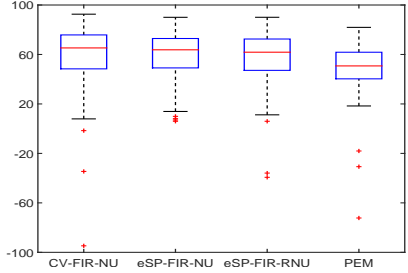}
\caption{Boxplot for Cram\'{e}r-Rao bound of 55\%.} 
\label{Boxplotl4} 
\end{center}
\end{figure}

The results of the experiment are provided in Table \ref{RNuclear}.
\begin{table}[h]
\centering
\caption{Average model fits of 150 randomly systems with white noise disturbance}
\begin{tabular}{|c|c|c|c|c|}
\hline
Cram\'{e}r Rao                  & CV-FIR-N & SPe-FIR-N & SPe-FIR-RN & PEM \\ \hline
\multicolumn{1}{|l|}{90\%} & 88.79 & 86.93        & 87.94 & 89.10 \\ \hline
\multicolumn{1}{|l|}{77\%} & 76.39 & 74.43        & 75.84 & 74.90 \\ \hline
\multicolumn{1}{|l|}{68\%} & 68.63 & 67.05        & 67.64 & 64.61 \\ \hline
\multicolumn{1}{|l|}{55\%} & 59.81 & 59.23        & 58.48 & 48.15 \\ \hline
\end{tabular}
\label{RNuclear}
\end{table}
\vspace{-0.1cm}
Boxplots of fit scores from the experiment with 4 different noise level are displayed in Figs. \ref{Boxplotl1} to \ref{Boxplotl4}. 
\vspace{-0.1cm}

\vspace{-0.1cm}

\vspace{-0.1cm}

From the above results, we can see that the SPe-FIR-N and SPe-FIR-RN methods show a competitive model fit when compared to the CV-FIR-N and PEM methods. In addition, the SPe-FIR-N shows a smaller number of bad estimates (outliers) than the CV-FIR-N and PEM methods. The SPe-FIR-RN shows an improvement on model fits with respect to the SPe-FIR-N. 

We also compute the average computational time of each estimation for the 4 methods. In this experiment, we run randomly 50 systems for 4 different noise levels, each system is estimated with one set of input and output data. The setting for the MATLAB  $\tt{fminsearch}$ command to search for the regularized parameter in the cross validation method is \texttt{optimset(`TolX',0.1,`TolFun',1e-4)}. The YALMIP toolbox by \cite{Yalmip:04} is used to implement all the regularization methods. This average computational time is presented in Table \ref{white_cotime}.

\begin{table}[H]
\centering
\caption{Average computational time of the 4 methods (in seconds)}

\begin{tabular}{|p{1.7cm}|p{1.7cm}|c|c|c|}
\hline
CV-FIR-N & SPe-FIR-N & SPe-FIR-RN & PEM \\ \hline
\multicolumn{1}{|l|}{251.12} & 12.88 & 145.88 & 1.14 \\ \hline
\end{tabular}
\label{white_cotime}
\end{table}

\vspace{-0.1cm}
Table \ref{white_cotime} shows that the computational time of the SPARSEVA-PEC regularization methods is much faster than the cross validation methods, which was expected. 

\vspace{-0.1cm}
\subsection{Coloured noise output disturbance}
\vspace{-0.2cm}

The experiment setting is same as in Sec. 5.1. The difference here is that in this experiment, we use a coloured noise output disturbance. The noise model has the same order as the system and is also generated by the $\tt{drss}$ command from Matlab with the same constraints on the poles.

In this case, we will use both FIR and ARX models to estimate the transfer function. The sample size in each run is $N=450$ and the number of parameters of the FIR model and ARX model is $n=35$ and $n_A=n_B=35$. 

The six methods that we compare are: 
\begin{enumerate}
\item
The nuclear norm regularization for FIR model with cross validation (CV-FIR-N).
\item
The nuclear norm regularization for ARX model with cross validation (CV-ARX-N).
\item
The nuclear norm regularization for FIR model with SPARSEVA-PEC (SPe-FIR-N).
\item
The reweighted nuclear norm regularization for FIR model with SPARSEVA-PEC (SPe-FIR-RN).
\item
The nuclear norm regularization for ARX model with SPARSEVA-PEC (SPe-ARX-N).
\item
The prediction error method (PEM) that uses the output-error model where it is assumed that the order of the system is already known.
\end{enumerate}
The results of the six numerical experiments are provided in Table \ref{RNuclear2}:

\begin{table}[h]
\centering
\caption{Average model fits of 150 randomly systems with coloured noise disturbance}
\begin{tabular}{|p{0.75cm}|p{0.75cm}|p{0.75cm}|p{0.75cm}|p{0.75cm}|p{0.75cm}|p{0.75cm}|}
\hline
Cram\'{e}r Rao                  & CV-FIR-N & CV-ARX-N & SPe-FIR-N & SPe-FIR-RN & SPe-ARX-N & PEM \\ \hline
\multicolumn{1}{|l|}{90\%} & 85.94 &  86.01       & 83.32 & 85.01 & 84.73 & 88.18 \\ \hline
\multicolumn{1}{|l|}{77\%} & 74.42 &   73.00      & 71.41 & 71.79 & 73.21 & 71.47 \\ \hline
\multicolumn{1}{|l|}{68\%} & 64.41 &   64.00      & 63.11 & 61.83 & 65.59 & 60.07 \\ \hline
\multicolumn{1}{|l|}{55\%} & 52.39 &    50.86     & 52.65 & 51.17 & 58.27 & 42.84 \\ \hline
\end{tabular}
\label{RNuclear2}
\end{table}

\begin{figure}[H]
\begin{center}
\includegraphics[width=66mm,height=45mm]{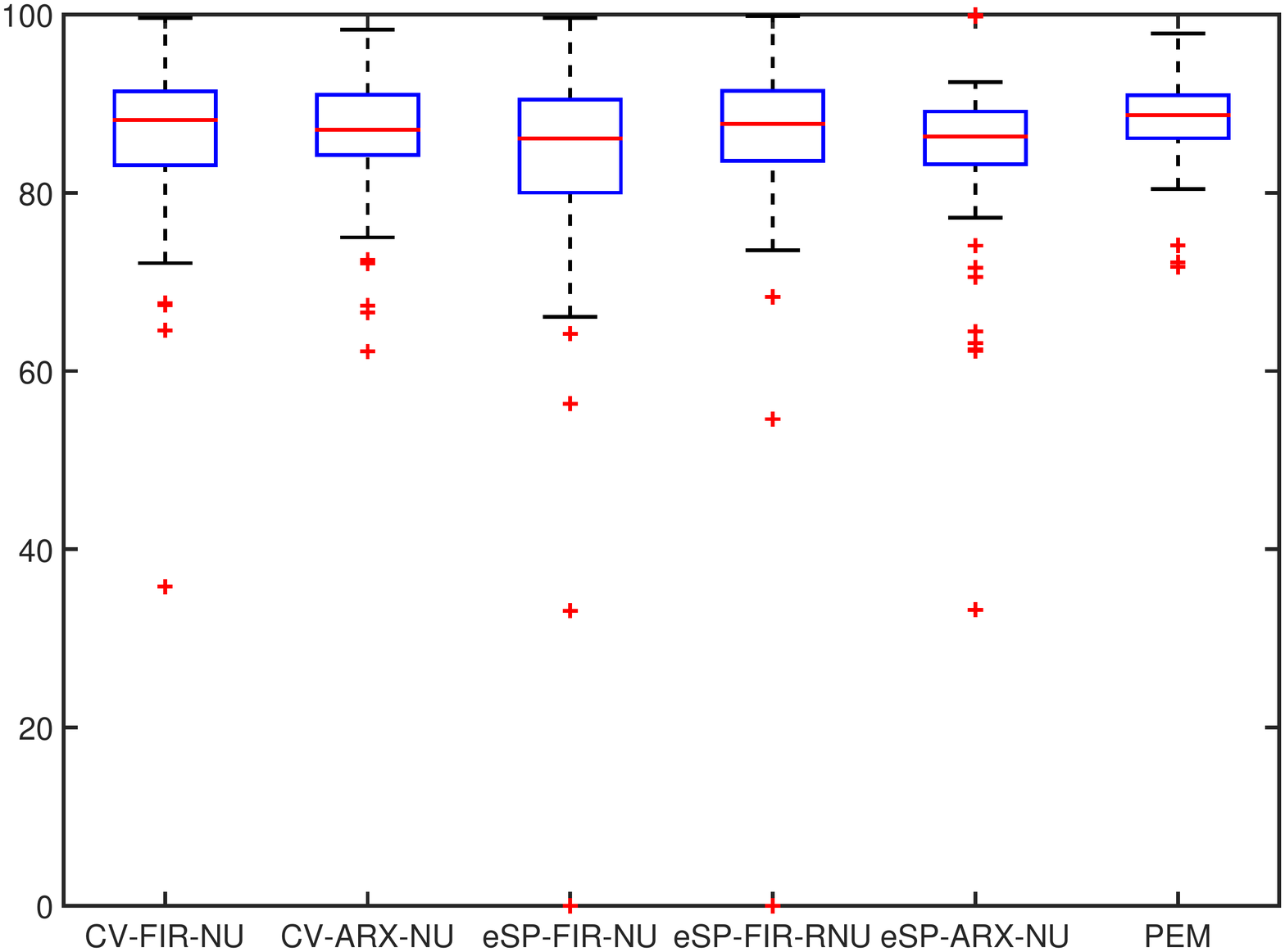}
\caption{Boxplot for Cram\'{e}r-Rao bound of 90\%.} 
\label{Boxplotl1c} 
\end{center}
\end{figure}

\begin{figure}[h]
\begin{center}
\includegraphics[width=66mm,height=45mm]{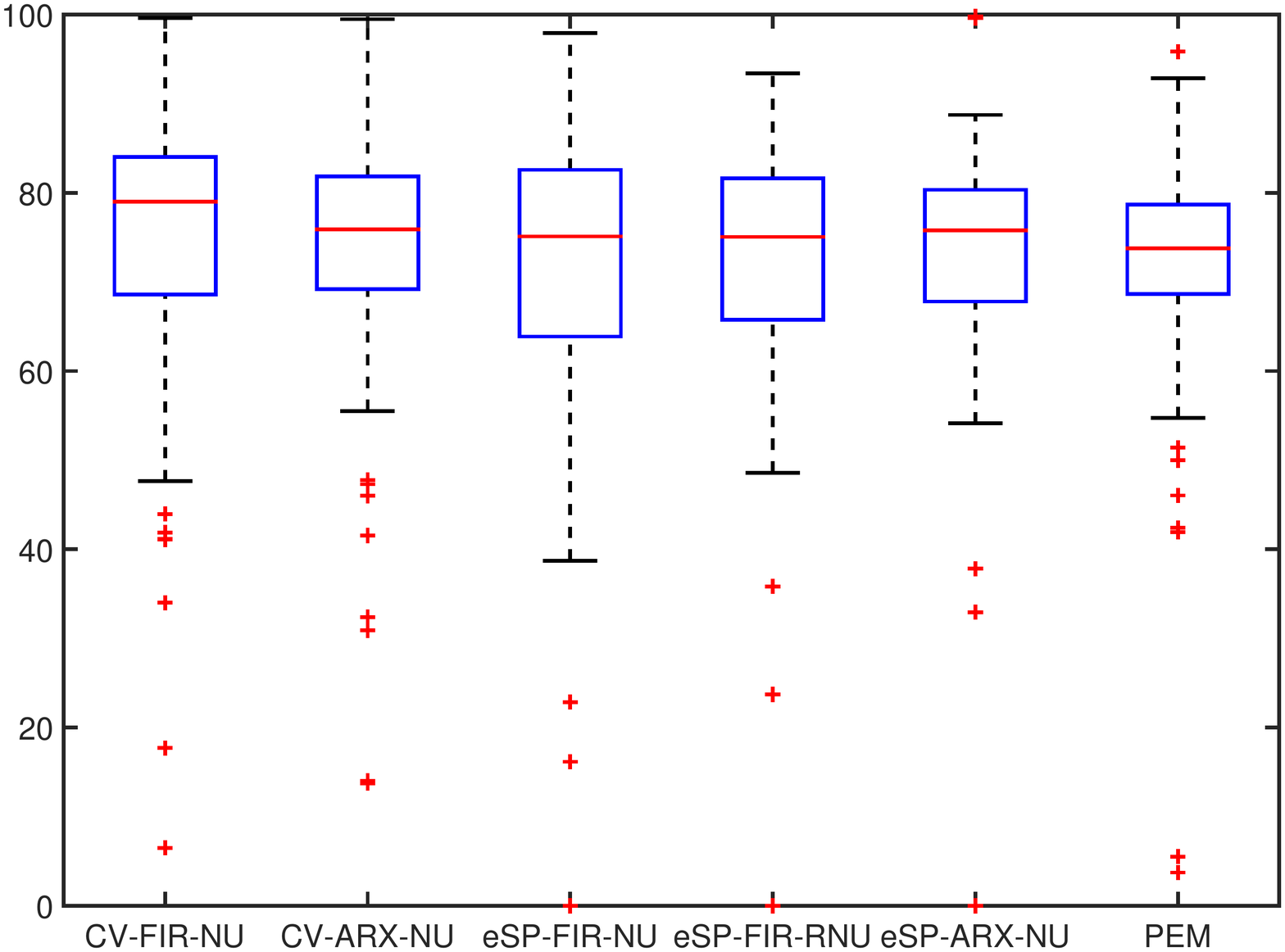}
\caption{Boxplot for Cram\'{e}r-Rao bound of 77\%.} 
\label{Boxplotl2c} 
\end{center}
\end{figure}

\begin{figure}[h]
\begin{center}
\includegraphics[width=66mm,height=45mm]{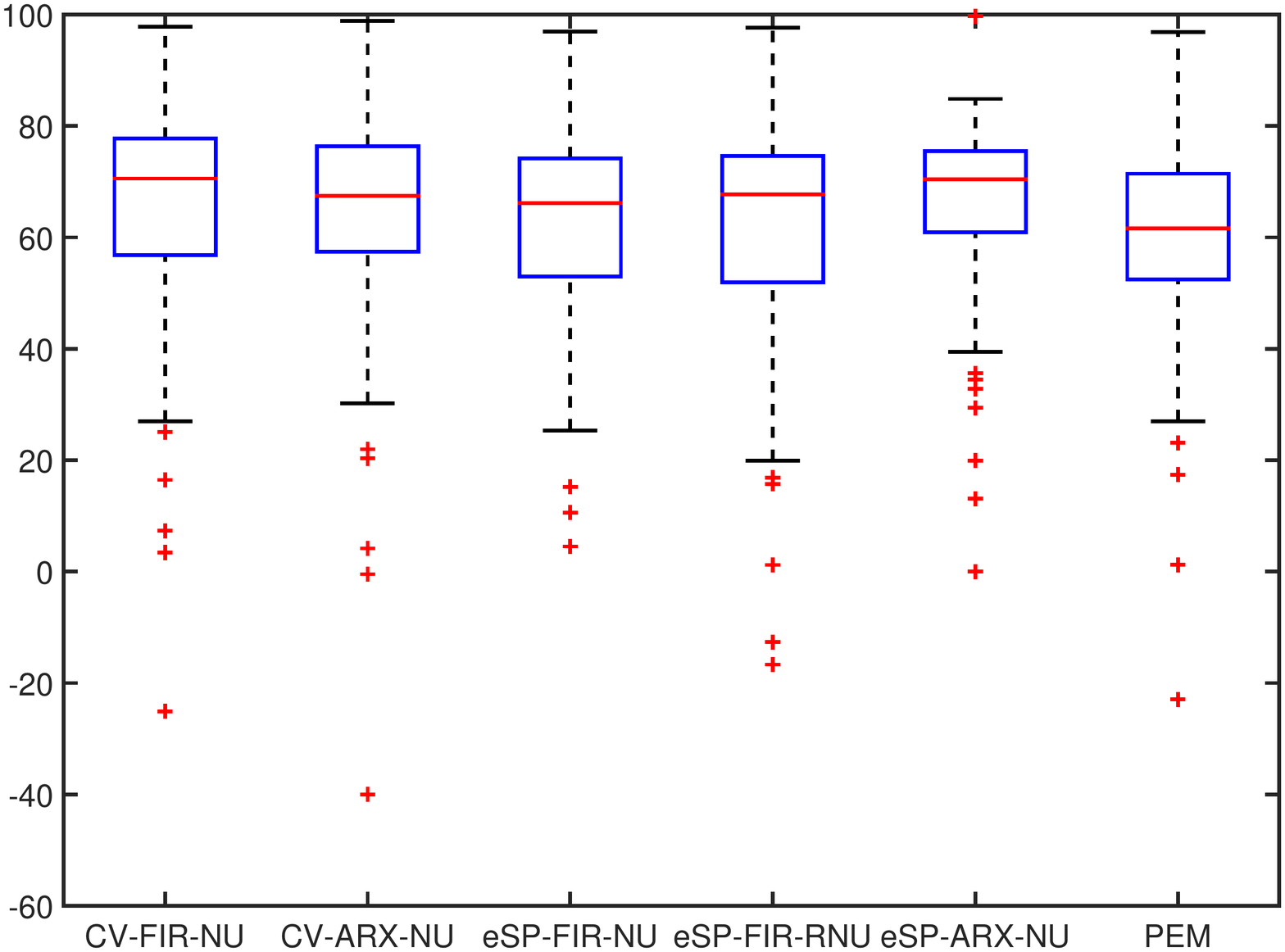}
\caption{Boxplot for Cram\'{e}r-Rao bound of 68\%.} 
\label{Boxplotl3c} 
\end{center}
\end{figure}

\begin{figure}[h]
\begin{center}
\includegraphics[width=66mm,height=45mm]{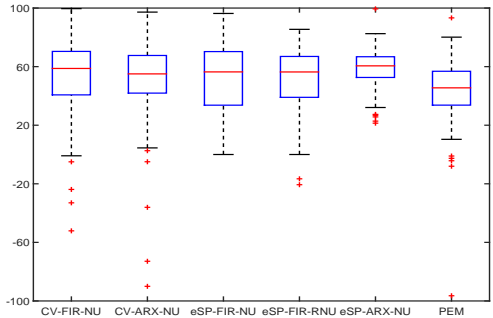}
\caption{Boxplot for Cram\'{e}r-Rao bound of 55\%.} 
\label{Boxplotl4c} 
\end{center}
\end{figure}

Boxplots of fit scores from the experiment with 4 different noise levels are displayed in Figs. \ref{Boxplotl1c} to \ref{Boxplotl4c}.

\vspace{-0.1cm}

\vspace{-0.1cm}

\vspace{-0.1cm}

\vspace{-0.1cm}

From the above results, we can see that the SPARSEVA-PEC methods with an FIR or an ARX model show a very good performance as compared to the cross validation methods and PEM. They especially perform better and have less bad estimates (outliers) when the SNR is low. Among all the SPARSEVA-PEC methods, the ARX structure gives better results than the FIR structure. 

The average computational time of each estimation for the 6 methods is also computed. The setting is the same as in the case of white noise disturbance. This average computational time is presented in Table \ref{color_cotime}.

\begin{table}[H]
\centering
\caption{Average computational time of the 6 methods (in seconds)}
\begin{tabular}{|p{0.75cm}|p{0.85cm}|p{0.75cm}|p{0.75cm}|p{0.75cm}|p{0.75cm}|p{0.75cm}|}
\hline
CV-FIR-N & CV-ARX-N & SPe-FIR-N & SPe-FIR-RN & SPe-ARX-N & PEM \\ \hline
\multicolumn{1}{|l|}{182.90} &   1152.91      & 13.19 & 142.89  & 27.79  & 0.83 \\ \hline
\end{tabular}
\label{color_cotime}
\end{table}
\vspace{-0.1cm}
The computational time of the SPARSEVA-PEC methods is again much faster than the cross validation methods.

\vspace{0.1cm}
\section{Conclusions}
\vspace{-0.1cm}

In this paper, a new method is developed to estimate a discrete time system using a high order FIR or ARX model with re-weighted nuclear norm regularization. The SPARSEVA framework is employed to avoid the computationally expensive searching of the tuning parameter in the re-weighted nuclear norm regularization problem. The paper also suggests the use of the PEC criterion for selecting the regularization parameter in the SPARSEVA framework. The numerical results show that the SPARSEVA-PEC nuclear norm regularization is competitive with respect to the accuracy of traditional technique of cross validation however the computational time is at least ten times faster.

\bibliography{ifacconf}             

\end{document}